\documentclass[english]{article}
\usepackage{amsxtra,amscd,amsthm,amsmath,amssymb}
\usepackage{longtable} \usepackage{bigstrut} \usepackage{enumerate}
\usepackage{cite}
\usepackage{verbatim}
\usepackage{graphicx}
\usepackage{xspace}

\usepackage{mathrsfs}

\def\vm{\vspace{-0.2cm}}
\def\vp{\vspace{0.1cm}}

\makeatletter
\def\setseccntfmt{\renewcommand{\@seccntformat}[1]{\S
    \csname the##1\endcsname.\hspace{1ex}}}
\def\setsubseccntfmt{\renewcommand{\@seccntformat}[1]{%
    (\csname the##1\endcsname)\hspace{0.5ex}}}
\renewcommand{\section}{\setseccntfmt\@startsection
  {section}{1}{0mm}{-\baselineskip}{0.5\baselineskip}{\bfseries\large}}
\def\presubsection{\setsubseccntfmt\@startsection
  {subsection}{2}{0mm}{-\baselineskip}{-0.5ex}{\bfseries\upshape}}
\def\tmpa{}\def\tmpb{}
\newcommand{\addspaceifnonempty}[1]{\def\tmpa{}\def\tmpb{#1}%
  \ifx\tmpa\tmpb{}\else{\hspace{0.5ex}#1\hspace{1ex}}\fi}
\renewcommand{\subsection}[1][]{\presubsection{\addspaceifnonempty{#1}}}
\def\presubsubsection{\setsubseccntfmt\@startsection
  {subsubsection}{3}{0mm}{-\baselineskip}{-0.5ex}{\bfseries\upshape}}
\renewcommand{\subsubsection}[1][]{\presubsubsection{\addspaceifnonempty{#1}}}

\usepackage{amsthm}

\theoremstyle{plain}
\numberwithin{equation}{subsection}

\newtheorem{ithm}[subsection]{Theorem}

\newtheorem{iDefinition}[subsection]{Definition}

\newtheorem{iconjecture}[subsection]{Conjecture}

\newcommand{\quash}[1]{}  
\newcommand{\nc}{\newcommand}
\nc{\on}{\operatorname}

\newcommand{\lps}{[\![}
\newcommand{\rps}{]\!]}
\newcommand{\llps}{(\!(}
\newcommand{\lrps}{)\!)}

\renewcommand{\phi}{\varphi}

\newcommand{\frakF}{{\mathfrak F}}

\newcommand{\eK}{{\mathbb K}}

\newcommand{\frakM}{{\mathfrak M}}

\newcommand{\frakQ}{{\mathfrak Q}}

\newcommand{\Sh}{{\rm Sh}}
\newcommand{\SSh}{{\mathscr S}}

\newcommand{\frakX}{{\mathfrak X}}

\newcommand{\DD}{{\mathbb D}}

\newcommand{\calB}{{\mathcal B}}

\newcommand{\calH}{{\mathcal H}}

\newcommand{\calP}{{\mathcal P}}

\newcommand{\VV}{{\mathbb V}}

\newcommand{\calZ}{{\mathcal Z}}

\newcommand{\ov}{\overline}

\nc{\al}{{\alpha}} \nc{\be}{{\beta}}

\nc{\ve}{{\varepsilon}} \nc{\Ga}{{\Gamma}}

\nc{\La}{{\Lambda}}

\def\0{\circ}

\newcommand{\BFf}{{\bar{\mathbb F}_p}}

\renewcommand{\cal}{\mathcal}

\newcommand{\E}{{\mathcal E}}
\newcommand{\A}{{\mathcal A}}

\renewcommand{\AA}{{\mathbb A}}

\newcommand{\GG}{{\mathbb G}}

\newcommand{\C}{{\mathbb C}}
\newcommand{\R}{{\mathbb R}}
\newcommand{\Q}{{\mathbb Q}}

\newcommand{\B}{{\mathcal B}}

\newcommand{\Ff}{{\mathbb F}}

\newcommand{\Gg}{{\mathcal G}}

\newcommand{\Gm}{{{\mathbb G}_{\rm m}}}
\newcommand{\Z}{{\mathbb Z}}
\newcommand{\F}{{\mathcal F}}
\newcommand{\ti}{\tilde}
\newcommand{\Spec}{{\rm Spec \, } }

 \renewcommand{\O}{{\mathcal O}}

\newcommand{\GL}{{\rm GL}}

\newcommand{\und}{\underline}

\def\thfill{\null\nobreak\hfill}

\def\endproof{\thfill\vbox{\hrule
  \hbox{\vrule\hbox to 5pt{\vbox to 5pt{\vfil}\hfil}\vrule}\hrule}}

\newcommand{\Om}{\Omega}

\newcommand{\GSp}{{\rm GSp}}

\newcommand{\fS}{{W_0\lps u\rps}}

\begin{document}
 
 \begin{center}
 {\Large \bf Arithmetic models for Shimura varieties}
 \end{center}
 \smallskip
 
  \begin{center}{\large \bf Georgios Pappas\footnote{Partially supported by  NSF grants DMS-1360733 and DMS-1701619. }
}
  \end{center}

\medskip
 

 \section*{ Introduction.}

Before the term became standard, certain   ``Shimura varieties" such as modular curves, Hilbert-Blumenthal varieties, and Siegel modular varieties, were already  playing an important role in number theory. Indeed, these are, respectively, quotients of the domains on which modular, Hilbert, and Siegel modular forms, are defined. In a series of groundbreaking works \cite{ShimuraAll, ShimuraBoundedI},
 Shimura initiated the arithmetic study of general quotients $\Gamma\backslash H$ of a hermitian symmetric domain $H$ by the action of a discrete congruence arithmetic group $\Gamma$ of holomorphic automorphisms of $H$. Such quotients are complex algebraic varieties and Shimura used the theory of moduli and of complex multiplication of abelian varieties to construct canonical models over explicit number fields for many of them.  

Deligne  reformulated and generalized   Shimura's theory and emphasized the group and motivic theoretic source of the constructions (\hspace{-0.001cm}\cite{DeligneTravauxShimura, DeligneCorvallis, HodgeCycles}).   
In Deligne's elegant definition, one starts with a pair $(\GG, X)$ of a connected reductive algebraic group $\GG$ defined over the rational numbers $\mathbb Q$ and a $\GG({\mathbb R})$-conjugacy class $X=\{h\}$ of an algebraic group homomorphism $h: {\mathbb S}\to \GG\otimes_{\mathbb Q}{\mathbb R}$. Here, ${\mathbb S}={\mathrm {Res}}_{\C/\R}(\Gm)$ is the algebraic torus over $\R$ whose real points are the group $\C^\times$. When the pair $(\GG, X)$ satisfies Deligne's conditions (\hspace{-0.001cm}\cite[(2.1.1.1)-(2.1.1.3)]{DeligneCorvallis}), we say that $(\GG, X)$ is a \emph{Shimura datum}. These conditions \quash{in particular,} imply that each connected component $X^+$ of $X$ is naturally a hermitian symmetric domain. 
The Shimura varieties 
$\Sh_\eK(\GG, X)$ are then the quotients
\[
\Sh_{\eK}(\GG, X)=\GG(\Q)\backslash X\times (\GG({\mathbb A}_f)/\eK)
\]
for $\eK$ an open compact subgroup of the finite adelic points $\GG({\mathbb A}_f)$ of $\GG$. Here, $\GG(\Q)$
acts diagonally on $X\times (\GG({\mathbb A}_f)/\eK)$, with the action on $X$  given by $\GG(\Q)\subset \GG(\R)$
followed by conjugation, and on $\GG({\mathbb A}_f)/\eK$ by $\GG(\Q)\subset \GG({\mathbb A}_f)$ followed by left translation.
Such a quotient is the disjoint union of a finite number of quotients of the form $\Gamma\backslash X^+$, where $\Gamma$ are discrete congruence arithmetic groups. Hence, $\Sh_\eK(\GG, X)$ has a natural complex analytic structure induced from that on $X$. In fact, by work of Baily and Borel, there is a quasi-projective complex algebraic structure on  $\Sh_\eK(\GG, X)$.  
Following the original work of 
Shimura and others, the existence of canonical models of $\Sh_K(\GG, X)$ over a number field 
(the ``reflex field" ${\mathbb E}={\mathbb E}(\GG, X)$) was shown in all generality by Borovoi \cite{BorovoiICM} and Milne (see \cite{MilneAnnArbor} and the references there). 

Many of the   applications of Shimura varieties in number theory depend on understanding 
 models of them over the ring of integers $\O_{\mathbb E}$ of the reflex field, or over localizations and completions of $\O_{\mathbb E}$. Indeed, perhaps the main application of Shimura varieties is to Langlands' program to associate Galois representations to automorphic representations. A related goal is to express the Hasse-Weil zeta functions of Shimura varieties as a product of automorphic $L$-functions.  
 
After some earlier work by Eichler, Shimura, Sato and Ihara, a general plan for realizing this goal was given by Langlands, first for the local factor of the  zeta function at a prime of good (\emph{i.e.} smooth) reduction.\footnote{At least for proper Shimura varieties;  in our brief report, we will ignore the issues arising from non-properness and the extensive body of work on compactifications.}  Langlands suggested expressing the numbers of points over finite fields of an integral model in terms of orbital integrals which appear in versions of the Arthur-Selberg trace formula. This was extended and realized in many cases, mainly by Kottwitz (\hspace{-0.001cm}\cite{KottShimuraAA,  KottJAMS}, see \cite{PicardBook}). It is now often referred to as the Langlands-Kottwitz method.
Langlands also considered the local factor of the zeta function for an example of a 
 Shimura surface at a prime of bad (non-smooth) reduction \cite{LanglandsBad}.
This example was treated carefully and the argument was extended to a larger class of Shimura varieties by Rapoport and Zink \cite{RapoportZinkMonodromy}
and by Rapoport \cite{RapoportAnnArbor}.  At primes of bad reduction, the singularities of the model have to be accounted for; the points need to be weighted by the semi-simple trace of Frobenius on the sheaf of nearby cycles. 
  In fact, the Galois action on  the nearby cycles can be used to   study the local Langlands correspondence, as in the work of Harris and Taylor \cite{HarrisTaylor}.  In relation to this,  Scholze  recently  extended the Langlands-Kottwitz method so that it can be applied, in principle at least, to the general case of bad reduction (\hspace{-0.001cm}\cite{ScholzeLKJAMS}).

A starting point of all the above is the existence of a reasonably well-behaved arithmetic model of the Shimura variety. For some time, such models could only be constructed for Shimura varieties of PEL type, \emph{i.e.} those given as moduli spaces of abelian varieties with additional polarization, endomorphism, and level structures, and over primes at which the level subgroup is ``parahoric'' (\emph{e.g.} \cite{RapZinkBook}). Recently, due  to advances in $p$-adic Hodge theory and in our understanding of the underlying group theory, the construction has been extended to most Shimura varieties of ``abelian type" (at good reduction 
by Kisin \cite{KisinJAMS}, see also earlier work of Vasiu \cite{VasiuAJM}; 
at general parahoric level in \cite{KisinPappas}). 
These Shimura varieties include most cases with $\GG$ a classical group.

The construction and properties of these   models is the subject of our report.
There are, of course, more uses for these  
in number theory besides in the Langlands program. For example, one could mention showing 
Gross-Zagier type formulas via intersection theory over the integers, or developments in the theory
of $p$-adic automorphic forms.
Here, we view their construction and study as a topic of its own and discuss it independently of applications. 

In fact, there are deep relations and analogies between this subject and the study of other
spaces of interest in number theory, representation theory, and the geometric Langlands program,
such as moduli spaces of bundles, or versions of affine Grassmannians and flag varieties.
Increasingly, these connections, especially with the geometric side of Langlands program, are taking center stage. 
For example, certain $p$-adically integral models of homogeneous spaces that appear as subschemes
of global (``Beilinson-Drinfeld") affine Grassmannians, the so-called ``local models", play an important role.
These local models also appear in the theory of deformations of Galois representations \cite{KisinFF}.
After first giving some background, we start by discussing local models. 
We then describe results on arithmetic models of Shimura varieties and their reductions, and finish with an account
of the local theory of Rapoport-Zink formal schemes. 
\smallskip

\emph{Acknowledgement:} I would like to thank my collaborators, especially M.~Kisin and X.~Zhu,  for making much of the progress 
reported on here possible, and M.~Rapoport for
generously sharing his knowledge over the years.

\section{ Recollections on $p$-adic groups.}\label{Local}

Let $G$ be a connected reductive group over the field of $p$-adic numbers $\Q_p$ for a prime number $p$.  Let $\bar\Q_p$ be an algebraic closure of $\Q_p$. We denote by $L$ the $p$-adic completion of the maximal unramified extension of $\Q_p$ in $\bar\Q_p$ and by $\O_L$ the integers of $L$. We will also denote by $k=\bar {\mathbb F}_p$ the algebraically closed residue field of $L$ and by $\sigma$ the automorphism of $L$ which lifts the Frobenius $x\mapsto x^p$.

\subsection{}  Recall that $g$, $g'\in G(L)$ are $\sigma$-conjugate, if there is $h\in G(L)$ with
$g'=h^{-1}g\sigma(h)$. We denote by $B(G)$ the set of $\sigma$-conjugacy classes of $G(L)$.
Recall the functorial surjective homomorphism
$
\kappa_G(L) : G(L)\to \pi_1(G)_I
$
of Kottwitz from \cite{KottIsocrystalsII}.
Here, ${\rm Gal}(\bar L/L)\simeq I\subset {\rm Gal}(\bar\Q_p/\Q_p)$ is the  inertia subgroup and $\pi_1(G)_I$ the inertia coinvariants of the algebraic fundamental group $\pi_1(G)$ of $G $ over $\bar\Q_p$. 
We denote the kernel of $\kappa_G(L)$ by $G(L)_1$.

\subsection{} Let $S$ be a maximal split torus of $G_L$. By Steinberg's theorem, $G_L$ is quasi-split and the centralizer of $S$
is a maximal torus $T$ of $G_L$. Denote by $N$ the normalizer of $T$. The quotient
$
\widetilde W=\widetilde W_{G, S}=N(L)/T(L)_1
$
is the \emph{Iwahori-Weyl group} associated to $S$. It is an extension of the relative Weyl group $W_0=N(L)/T(L)$
by $\pi_1(T)_I$. Since $\pi_1(T)=X_*(T)$ (the group of cocharacters of $T$ over $\bar L$), we obtain
an exact sequence
$
1\to X_*(T)_I\to \widetilde W\to W_0\to 1.
$

\subsection{} 
Suppose that $\{\mu\}$ is the conjugacy class of a cocharacter $\mu: \Gm_{\bar\Q_p}\to G_{\bar\Q_p}$.
Then $\{\mu\}$ is defined over the \emph{local reflex field} $E$ which is a finite extension of $\Q_p$ contained in $\bar E=\bar\Q_p$. Denote by $\O_E$ the integers of $E$ and by $k_E$ its residue field. There is a corresponding 
  homogeneous space $G_{\bar E}/P_{\mu^{-1} }$
 which   has a canonical model $X_\mu=X(\{\mu\})$ defined over $E$. Here, $P_\nu$ denotes the parabolic subgroup  
 that corresponds to the cocharacter $\nu$.

A pair $(G, \{\mu\})$ of a connected reductive group $G$ over $\Q_p$, together with a conjugacy class $\{\mu\}$ of a cocharacter $\mu$ as above, is a \emph{local Shimura pair},\footnote{Compare this with the term ``local Shimura datum" used in \cite{RapoportViehmann}. A local Shimura datum $(G, \{\mu\}, [b])$ also includes the choice of a $\sigma$-conjugacy class $[b]$.} if $\mu$ is \emph{minuscule} (\emph{i.e.},
for any root $\alpha$ of $G_{\bar\Q_p}$, $\langle \alpha, \mu\rangle\in \{-1,0,1\}$.)

\subsection{} 
We  denote by $\calB(G, \Q_p)$ the (extended) Bruhat-Tits building of $G(\Q_p)$
\cite{BTI, BTII, TitsCorvallis}. The group $G(\Q_p)$ acts on $\calB(G, \Q_p)$ on the left.
If $\Omega $ is a   subset of $ \cal B(G, \Q_p)$,
we   write 
$
G(\Q_p)_\Omega=\{g\in G(\Q_p)\ |\ g\cdot y=y, \hbox{\rm\ for all\ } y\in \Omega\}
$
for the pointwise stabilizer  of $\Omega$.  
Similarly,  
we have the subgroup $G(L)_\Om$
of $G(L)$. 

By the main result of \cite{BTII}, if  $\Omega$ is bounded and contained in an apartment, there is a smooth affine group scheme $\Gg_\Om$ over $\Spec(\Z_p)$
with generic fiber $G$ and with $\Gg_\Om(\O_L)=G(L)_\Om$, which is uniquely characterized by these properties. 
By definition, the ``connected stabilizer" $G(\Q_p)_\Omega^\circ$ is $\Gg^\circ_\Om(\Z_p),$ where $\Gg^\circ_\Om$ is the connected component of $\Gg_\Om.$ It is a subgroup of finite index in $G(\Q_p)_\Omega$.
When $\Omega=\{x\}$ is a point, we simply write $G(\Q_p)_x$ and $G(\Q_p)_x^\circ$.
If $\Omega$ is an open facet and $x\in \Omega$, then $G(\Q_p)_\Omega^\circ=G(\Q_p)_x^\circ$.

A \emph{parahoric} subgroup of $G(\Q_p)$ is any subgroup which is the connected stabilizer $G(\Q_p)_x^\circ$ of some point $x$ 
in $\B(G, \Q_p)$ as above. 
These are open compact subgroups of $G(\Q_p)$.
We call $\Gg^\circ_x$ a ``parahoric group scheme".

We now recall some more terms and useful facts  (\hspace{-0.001cm}\cite{TitsCorvallis}, \cite{HainesRapoportAppendix}):
A point $x\in \B(G, \Q_p)$ is \emph{hyperspecial}, when $\Gg_x$ is reductive; then $\Gg_x=\Gg_x^\circ$.
Hyperspecial points exist if and only if  
$G$ is \emph{unramified} over $\Q_p$, \emph{i.e.} if $G$ is quasi-split and splits over an unramified extension of $\Q_p$.
When $x$ is a point in an open alcove, then the parahoric $G(\Q_p)_x^\circ$ is called an \emph{Iwahori} subgroup.
Iwahori subgroups exist for all $G$.  
Each maximal open compact subgroup of $G(\Q_p)$ contains a parahoric subgroup with finite index. 
If the group  $\pi_1(G)_I$ 
has no torsion, then $G(\Q_p)_x=G(\Q_p)^\circ_x$, for all $x$. Then, all maximal open compact subgroups are parahoric.

\subsection{}
Let $K=G(\Q_p)_x^\circ $ be parahoric 
and denote by $\ti K$ the corresponding
 parahoric subgroup $G(L)_x^\circ$ of $G(L)$; we have $\ti K=\Gg_x^\circ(\O_L)$.
 
 Suppose $x$ lies in the apartment associated to $S$.
A choice of an alcove $C\subset \B(G, L)$ contained in that apartment
provides the Iwahori-Weyl group $\widetilde W=\widetilde W_{G,S}$ with a Bruhat partial order $\leq $.
Suppose $x\in C$ and set
$
\widetilde W^{\ti K} = (N(L)\cap  \ti K)/T(L)_1, 
$
which is a subgroup of $\widetilde W$.  
The inclusion $N(L)\subset G(L)$  induces \cite{HainesRapoportAppendix} a bijection  
\[\vspace{-0.15cm}
  \widetilde W^{\ti K}\backslash \widetilde W/\widetilde W^{\ti K}\xrightarrow{\sim} \ti K\backslash G(L)/\ti K.
\] 
There is also a partial order $\leq $ on these double cosets: Set  $ [w]=\widetilde W^{\ti K} w\widetilde W^{\ti K}$. 
Then  $[w_1]\leq [w_2]$  
 if and only if there are $w_1'\in [w_1]$, $w'_2\in [w_2]$, with $w'_1\leq w_2'$.
 
 We refer the reader to \cite{PRS} for the definition of the \emph{$\{\mu\}$-admissible subset}
 \[
 {\rm Adm}_{\ti K}(\{\mu\}) \subset \widetilde W^{\ti K}\backslash \widetilde W/\widetilde W^{\ti K}\simeq \ti K\backslash G(L)/\ti K
 \]
  of Kottwitz and Rapoport. This is a finite set which has the following property: If $[w]\in {\rm Adm}_{\ti K}(\{\mu\})$ and $[w']\leq [w]$, then
 $[w']\in {\rm Adm}_{\ti K}(\{\mu\})$.

\subsection{}\label{ADL}

Let us continue with the above set-up. 
   The   \emph{affine Deligne-Lusztig  set   $X_K(\{\mu\}, b)$}
   associated to $G$, $\{\mu\}$, $K$,
and $b\in G(L)$, is the subset of $G(L)/\ti K$ consisting of those cosets $g\ti K$ for which
$g^{-1}b\sigma(g)\in \ti Kw\ti K$, for some $[w]\in {\rm Adm}_{\ti K}(\{\mu\})$.

If $b'$ and $b$ are $\sigma$-conjugate $b'=h^{-1}b\sigma(h)$, then $g\ti K\mapsto hg\ti K$ gives a bijection
$
 X_K(\{\mu\}, b)\xrightarrow{\sim} X_K(\{\mu\}, b').
$
 The group $J_b(\Q_p)=\{j\in G(L)\ |\ j^{-1}b\sigma(j)=b\}$ acts on $X_K(\{\mu\}, b)$ by $j\cdot g\ti K=jg\ti K$. 
 Set $f=[k_E:\Ff_p]$. The identity $\Phi_E(g\ti K)=b\sigma(b)\cdots \sigma^{f-1}(b)\sigma^f(g)\ti K$
 defines a map $\Phi_E: X_K(\{\mu\}, b)\to X_K(\{\mu\}, b)$.
 
\subsection{}\label{aff} 
Our last reminder is of the affine Grassmannian of $G$. By definition, the \emph{affine Grassmanian} ${\rm Aff}_{G}$ of $G$  is the  ind-projective ind-scheme 
 over $\Spec(\Q_p)$ which represents (\emph{e.g.} \cite{PappasRaTwisted}) the fpqc sheaf associated to the quotient presheaf given by $R\mapsto G(R\llps t\lrps)/G(R\lps t \rps)$.

\section{\bf Local models.}

Let $(G,  \{\mu\})$ be a local Shimura pair and $K $ a parahoric subgroup 
of $G(\Q_p)$.  
Denote by $\Gg$ the corresponding parahoric group scheme  
 over $\Spec(\Z_p)$ so that $K=\Gg(\Z_p)$ and set again $\ti K=\Gg(\O_L)$. A form of the following appears 
 in \cite{RapoportGuide}.
 
 \begin{iconjecture}\label{conjLocal} There exists a projective and flat scheme ${\rm M}_K({G, \{\mu\}})$ over $\Spec(\O_E)$
 which supports an action of $\Gg\otimes_{\Z_p}\O_E$ and is such that:
 
 (i) The generic fiber ${\rm M}_K({G, \{\mu\}})\otimes_{\O_E} E$ is $G_E$-isomorphic to $X_\mu$,
 
 (ii) There is a $\ti K$-equivariant bijection between $ {\rm M}_K({G, \{\mu\}})(k)$ and
 \vspace{-0.1cm}
 \[
 \{g\tilde K\in G(L)/\ti K\ |\ \ti K g\ti K\in {\rm Adm}_{\ti K}(\{\mu\})\}\subset G(L)/\ti K.
 \]

\end{iconjecture}

The scheme ${\rm M}_K({G, \{\mu\}})$ is a \emph{local model} associated to $(G,\{\mu\})$ and $K$.

\subsection{}
We now consider the following ``tameness" condition:
 \vp
 
\noindent \emph{ (T) The group $G$ splits over a tamely ramified extension of $\Q_p$ and $p$ does not divide the order of the fundamental group ${\rm \pi}_1(G^{\rm der})$ of the derived group.  }
\vp

Under the assumption (T), Conjecture \ref{conjLocal} is shown in \cite{PaZhu}. The construction 
of the local models in \cite{PaZhu} is as follows. 
\smallskip
 
\noindent (I)  We first construct a smooth affine group scheme $\und \Gg$ over $\Z_p[u]$ which, among other properties, satisfies: 
 
 1)  the restriction of $\und\Gg$ over $\Z_p[u,u^{-1}]$ is reductive, 
 
 2) the base change $\und\Gg\otimes_{\Z_p[u]}\Z_p$, given by $u\mapsto p$, is isomorphic to the parahoric group scheme $\Gg$.

 For this, we use that $\Z_p[u^{\pm 1}]=\Z_p[u,u^{-1}]\to \Q_p$, given by $u\mapsto p$,  identifies the \'etale fundamental group of $\Z_p[u^{\pm 1}]$ with the tame quotient of the Galois group of $\Q_p$; our tameness hypothesis enters this way. We first obtain the restriction  $\und \Gg|_{\Z_p[u^{\pm 1}]}$ from its constant Chevalley form $H\otimes_{\Z_p}\Z_p[u^{\pm 1}]$
 by using, roughly speaking, the ``same twist" that gives $G$ from its constant form $H\otimes_{\Z_p}\Q_p$.
 This  is quite straightforward when $G$ is quasi-split and the twist is given using a diagram automorphism; the general case uses explicit Azumaya algebras over $\Z_p[u^{\pm 1}]$. The extension $\und\Gg$ to $\Z_p[u]$ of the reductive group scheme
 $\und \Gg|_{\Z_p[u^{\pm 1}]}$ is then given by generalizing some of the constructions of Bruhat and Tits in \cite{BTII} to this two-dimensional set-up.
 
  \smallskip
  
\noindent (II)   Consider the functor that sends $\phi:\Z_p[u]\to R$ to the set of isomorphism classes of pairs $(\E, \beta)$, 
with $\E$ a $\und\Gg$-torsor over $\AA^1_R$ and $\beta$ a section of $\E$ over $\Spec(R[u, (u-\phi(u))^{-1}])$. 
We show that this is represented by an ind-projective ind-scheme 
(the \emph{global} affine Grassmannian, \emph{cf.} \cite{BeDr})
${\rm Aff}_{\und \Gg, \AA^1}\to
\AA^1=\Spec(\Z_p[u])$.  
Using the construction of $\und\Gg$ and \cite{BeauvilleLaszlo}, we obtain an equivariant isomorphism
$$
 {\rm Aff}_{G}\xrightarrow{\sim}{\rm Aff}_{\und \Gg, \AA^1}\otimes_{\Z_p[u]} \Q_p
$$
where the base change is  given by $u\mapsto p$.
(The above isomorphism is induced by $t\mapsto u-p$, with the notation of  \S \ref{aff} for ${\rm Aff}_{G}$.)
 \smallskip

 \noindent (III) \ Since $\Gm=\Spec(\Q_p[t, t^{-1}])$, our cocharacter $\mu$ defines $\mu(t)\in G(\bar \Q_p\llps t\lrps)$. In turn, this gives a $\bar \Q_p$-point  $[\mu(t)]=\mu(t)G(\bar \Q_p\lps t \rps)$ of   ${\rm Aff}_G$.
 
  Since $\mu$ is minuscule and the conjugacy class $\{\mu\}$ is defined over $E$, the orbit
 \[
G(\bar \Q_p\lps t\rps )\cdot [\mu(t)]\subset G(\bar \Q_p\llps t\lrps)/G(\bar \Q_p\lps t \rps)={\rm Aff}_{G}(\bar \Q_p),
\]
 is equal to the set of $\bar \Q_p $-points of a
\emph{closed} subvariety $S_\mu$ of  ${\rm Aff}_{G, E}:={\rm Aff}_{G}\otimes_{\Q_p}E$. 
The stabilizer $H_{\mu^{-1}}:=G(\bar \Q_p\lps t\rps ) \cap \mu(t)G(\bar \Q_p\lps t\rps ) \mu(t)^{-1}$ of $[\mu(t)]$
is the inverse image of $P_{\mu^{-1}}(\bar \Q_p)$ under  $G(\bar \Q_p\lps t\rps ) \to G(\bar \Q_p)$
given by reduction modulo $t$. 
This gives a $G_E$-equivariant isomorphism 
$
X_\mu\xrightarrow{\sim} S_\mu 
$
which we can use to identify $X_\mu$ with $S_\mu\subset {\rm Aff}_{G, E}$.

\begin{iDefinition}\label{defLocal}
We define 
$
{\rm M}_K({G, \{\mu\}})
$
to be the (flat, projective) scheme over $\Spec(\O_E)$ given by the reduced Zariski closure of the image of
\[
X_\mu\xrightarrow{\sim} S_\mu\subset {\rm Aff}_{G, E}\xrightarrow{\sim} {\rm Aff}_{\und \Gg, \AA^1}\otimes_{\Z_p[u]} E
\]
in the ind-scheme ${\rm Aff}_{\und \Gg, \AA^1}\otimes_{\Z_p[u]} \O_E$.
\end{iDefinition}

\subsection{}\label{resultsLocal}

The next includes the main general facts about the structure of ${\rm M}_K({G, \{\mu\}})$
and can be extracted from the results of \cite{PaZhu}. 

\begin{ithm}\label{normalThm}
Suppose that $(T)$ holds. The scheme ${\rm M}_K({G, \{\mu\}})$ of Definition \ref{defLocal}
satisfies Conjecture \ref{conjLocal}. In addition:

 a) The scheme
${\rm M}_K({G, \{\mu\}})$
is normal.

b) The geometric special fiber ${\rm M}_K({G, \{\mu\}})\otimes_{\O_E} k$ is reduced and admits a $\Gg\otimes_{\Z_p}k$-orbit stratification by locally closed and smooth  strata $S_{[w]}$, 
with 
$
S_{[w]}(k)\simeq\{g\ti K\in G(L)/\ti K\ |\ [w]=\ti Kg\ti K\},
$
for each $[w]\in  {\rm Adm}_{\ti K}(\{\mu\})$.

  c)  The closure $\ov S_{[w]}$ of each stratum is normal and Cohen-Macaulay
and equal to the union
$
\bigcup_{[w']\leq [w]} S_{[w']},
$
where $\leq $ is given by the Bruhat order.

\end{ithm}

One main ingredient is the proof, by  Zhu \cite{ZhuCoh}, of the coherence conjecture of the author and 
 Rapoport \cite{PappasRaTwisted}. The coherence conjecture is a certain numerical
 equality in the representation theory of (twisted) Kac-Moody groups. Its statement is, roughly speaking, 
 independent of the characteristic, and so enough to show in the function field case where more tools are available.

Before \cite{PaZhu}, there have been various, often ad hoc,   constructions of 
local models 
obtained by using
variants of linked Grassmannians (\hspace{-0.001cm}\cite{RapZinkBook}, \cite{PaJAG}; see the survey \cite{PRS} and the references within, and \S \ref{exampleSympl2} below
for an example). 

An extension of the above construction of local models, to the case
$G$ is the restriction of scalars of a tame group from a wild extension, is given by   Levin in \cite{LevinLM}. One would expect that 
a completely general construction of the local model ${\rm M}_K(G, \{\mu\})$ can be given by using a hypothetical version of the affine Grassmannian in which $R\mapsto \Gg(R\lps u\rps)$ is replaced by the  functor $R\mapsto \Gg(W(R))$. Here $W(R)$ is the ring of $p$-typical Witt vectors of $R$. (See \S \ref{displays}
for a related construction when $\Gg$ is reductive.) Such Witt affine Grassmannians are defined 
 in characteristic $p$  by Zhu \cite{ZhuWitt} (but only ``after perfection") and further studied by Bhatt and Scholze \cite{BhattScholzeWitt}; other variations have been used by Scholze.

 \section{\bf  Global theory: Arithmetic models.}

Let $(\GG, \{h\})$ be a Shimura datum as in the introduction. Define $\mu_h: \Gm_\C\to \GG_\C$
by $\mu_h(z)=h_\C(z, 1)$. The $G(\C)$-conjugacy class $\{\mu_h\}$ is defined over the 
reflex field $\mathbb E\subset \bar\Q\subset \C$. Write ${\mathbb Z}_s$ for the maximal subtorus of the center $Z(\GG)$  
which is $\mathbb R$-split but which has no $\Q$-split subtorus  and set $\GG^c=\GG/{\mathbb Z}_s$.

\subsection{}\label{sympl}
Consider a $\Q$-vector space $\mathbb V$   with a perfect alternating bilinear pairing $\psi: \mathbb V\times \mathbb V\to \Q$. Set
$\dim_\Q(\mathbb V)=2g$. Let $\GG={\rm GSp}(\mathbb V, \psi)$ be the group of symplectic similitudes and 
$X=S^{\pm}$ the Siegel double space. This is the set of homomorphisms $h: \mathbb S\to \GG_{\R}$
such that: (1) The $\C^\times$-action given by $h(\R): \C^*\to \GG(\R)$ gives on $\mathbb V_\R$ a Hodge structure of type $\{(-1,0), (0, -1)\}$:
$
\mathbb V_\C=\mathbb V^{-1. 0}\oplus \mathbb V^{0, -1},
$
and, (2) the form $(x, y)\mapsto \psi(x, h(i)y)$ on $\mathbb V_{\R}$ is positive (or negative) definite. 
The pairs $({\rm GSp}(\mathbb V, \psi), S^\pm)$ give the most important examples 
of Shimura data. By Riemann's theorem, if $\VV_\Z$ is a $\Z$-lattice in $\mathbb V$, and $h\in S^{\pm }$, the quotient torus $\mathbb V^{-1, 0}/\VV_\Z$
is a complex abelian variety of dimension $g$. This leads to an interpretation of the   Shimura varieties
$\Sh_{\eK}({\rm GSp}(\mathbb V, \psi), S^\pm)$ as moduli spaces of (polarized) abelian varieties of dimension $g$
with level structure.

\subsection{}\label{Hodgedef}
A Shimura datum $(\GG, X)$ with $X=\{h\}$ is of \emph{Hodge type}, when there 
is a symplectic space $(\mathbb V,\psi)$ over $\Q$ and a closed embedding 
$
\rho: \GG\hookrightarrow {\rm GSp}(\mathbb V, \psi)
$
such that the composition $\rho\circ h$ lies in the Siegel double space $S^{\pm }$. Then, we also have
$\GG=\GG^c$. In this case, the Shimura varieties $\Sh_\eK(\GG, X)$
 parametrize abelian varieties together with (absolute) Hodge cycles (see below). 
 A special class of Shimura data of Hodge type are those of \emph{PEL type} \cite[\S 5]{KottJAMS}, \cite[Ch. 6]{RapZinkBook}. 
For those, the corresponding Shimura varieties $\Sh_{\eK}(\GG, X)$ are (essentially) moduli schemes of abelian varieties together with polarization, endomorphisms and level structure.

A Shimura datum $(\GG, X)$ is of \emph{abelian type} if there is a datum of Hodge type 
$(\GG_1, X_1)$ and a central isogeny $\GG^{\rm der}_1\to \GG^{\rm der}$ between derived groups
which induces an isomorphism $(\GG^{\rm ad}_1, X^{\rm ad}_1)\xrightarrow{\sim} (\GG^{\rm ad}, X^{\rm ad})$.
Here, the superscript $^{\rm ad}$ denotes passing to the adjoint group: If $X=\{h\}$, then $X^{\rm ad}=\{h^{\rm ad}\}$
with $h^{\rm ad}$   the composition of $h:\mathbb S\to \GG_\R$ with $\GG_\R\to \GG^{\rm ad}_\R$.
Roughly speaking, most Shimura data $(\GG, X)$ with $\GG$ a classical group are of abelian type.

\subsection{} 
Fix a prime number $p$ and a prime $\mathfrak p$   of the reflex field $\mathbb E$ above $p$ which is obtained from an embedding $\bar\Q\hookrightarrow \bar\Q_p$. Let $E$ be the completion of $\mathbb E$ at $\mathfrak p$ and set $G=\GG_{\Q_p}$. Then $\{\mu_h\}$ gives a conjugacy class $\{\mu\}$ of $G$ defined over the local reflex field $E$ and $(G, \{\mu\})$ is a local Shimura pair.

Fix  a \emph{parahoric} subgroup $K=K_p\subset \GG(\Q_p)=G(\Q_p)$. For any open compact subgroup $K^p\subset \GG(\AA_f^p)$ we can consider the Shimura variety $\Sh_{\eK}(\GG, X)$ over $\mathbb E$, where $\eK=K_pK^p\subset \GG(\AA_f)$.  

We now assume that a local model ${\rm M}_K(G, \{\mu\})$ as in Conjecture \ref{conjLocal}
is given for $(G, \{\mu\})$ and 
parahoric $K$.

\begin{iconjecture}\label{conjintegral}
There is a scheme $\SSh_{K_p}(\GG, X)$ over $\Spec(\O_E)$
which supports a right action of $\GG(\AA^p_f)$ and has the following properties:

a) Any sufficiently small open compact subgroup $K^p\subset \GG(\AA_f^p)$ acts freely on 
$\SSh_{K_p}(\GG, X)$, and the quotient $\SSh_{\eK}(\GG, X):=\SSh_{K_p}(\GG, X)/K^p$ is a scheme of finite type over $\O_E$ which extends $\Sh_{\eK}(\GG, X)\otimes_{\mathbb E}E$. We have
$
\SSh_{K_p}(\GG, X)=\varprojlim\nolimits_{K^p} \SSh_{K_pK^p}(\GG, X)
$
where the limit is over all such $K^p\subset \GG(\AA_f^p)$.

b) For any discrete valuation ring $R\supset \O_E$ of mixed characteristic $(0, p)$, the map $\SSh_{K_p}(\GG, X)(R)\to \SSh_{K_p}(\GG, X)(R[1/p])$ is a bijection.

c) There is a smooth morphism of stacks over $\Spec(\O_E)$
\[
\lambda: \SSh_{K_p}(\GG, X)\to [(\Gg^c\otimes_{\Z_p}\O_E)\backslash {\rm M}_K(G, \{\mu\})]
\]
which is invariant for the $\GG(\AA^p_f)$-action on the source and is such that the base change $\lambda_E$
is given by the canonical principal $\GG^c_E$-bundle over $\Sh_{K_p}(\GG, X)\otimes_{\mathbb E}E$
(\hspace{-0.001cm}\cite[III, \S 3]{MilneAnnArbor}). Here, we set  $\Gg^c=\Gg/\calZ_s$,  where $\calZ_s$ is the Zariski closure of the central torus $Z_s\subset G$ in $\Gg$.  
\end{iconjecture}

It is important to record that the existence of the smooth $\lambda$ as in (c) implies:

\smallskip

\noindent \emph{(c') For each closed point $x\in \SSh_{K_p}(\GG, X)$, there is a closed point $y\in{\rm M}_K(G, \{\mu\})$
and \'etale neighborhoods $U_x\to \SSh_{K_p}(\GG, X)$ of $x$ and $V_{y}\to {\rm M}_K(G, \{\mu\})$ of $y$, 
which are
isomorphic over $\O_E$. }

\smallskip
 
The significance of (c') is that the singularities of ${\rm M}_K(G, \{\mu\})$ control the singularities of $\SSh_{K_p}(\GG, X)$, and so, by (a), also of the integral models
$\SSh_{\eK}(\GG, X)$ for $\eK=K_pK^p$, with $K^p$ sufficiently small. In what follows, we will
often assume that $K^p$ is sufficiently small without explicitly saying so.

The ``extension property'' (b) ensures that the special fiber of $\SSh_{K_p}(\GG, X)$ is sufficiently
large, in particular, it cannot be empty. Unfortunately, it is not clear that properties (a)-(c)
uniquely characterize the schemes $\SSh_{K_p}(\GG, X)$. We will return to this question later.

\subsection{}\label{exampleSympl2} 
Before we give some general results we discuss the illustrative example of the Siegel Shimura datum $({\rm GSp}(\VV,\psi), S^\pm)$ as in \S\ref{sympl}.

Fix a prime $p$. Set $V=\mathbb V_{\Q_p}$ and denote the induced form on $V$ also by $\psi$. For a $\Z_p$-lattice $\La$ in $V$, its dual is
$
\La^\vee=\{ x\in  V\ |\ \psi(x, y)\in \Z_p\, , \hbox{\rm {for all $y\in \La$}}\}.
$
Now choose a chain of lattices $\{\La_{i}\}_{i\in \Z}$ in $V$ with the following properties: $\La_i\subset \La_{i+1}$ and $\La_{i-m}=p\La_i$, for all $i$ and a fixed ``period" $m>0$, 
 and 
$\La^\vee_i=\La_{-i+a}$, for all $i$ and (fixed) $a=0$ or $1$. 
  The chain $\{\Lambda_i\}_{i\in \Z}$ gives rise to a point in the building of $G={\rm GSp}(V, \psi)$ over $\Q_p$. The corresponding parahoric $K$ is the subgroup of all $g\in {\rm GSp}(V, \psi)$ such that $g\Lambda_i=\Lambda_i$, for all $i$. Every parahoric subgroup of $G(\Q_p)$ is obtained from such a lattice chain.
 
 In this case, schemes $\SSh_{K}({\rm GSp}(\mathbb V, \psi), S^{\pm})$ that satisfy the conjecture
 are given by moduli spaces as we will now explain. 
Recall that the category $AV_{(p)}$ of abelian schemes up to prime-to-$p$ isogeny has objects abelian schemes
and morphisms given by tensoring the Hom groups by $\Z_{(p)}$.
We consider the functor which associates to a $\Z_{p}$-algebra $R$ the set of isomorphism classes of $m$-tuples $\{(A_i, \alpha_i, \ti \nu, \eta)\}_{i\in \Z/m\Z}$ where (see \cite[Ch. 6]{RapZinkBook} for details):

1) $A_i$ are 
 up to prime-to-$p$ isogeny projective abelian schemes  over $R$, 
 
2)  $\alpha_i: A_i\to A_{i+1}$ are isogenies of height $\log_p([\La_{i+1}:\La_i])$ in $AV_{(p)}$
 such that the compositions $\alpha_{r+m-1}\cdots \alpha_{r+1}\alpha_r=p$, for all $0\leq r< m$,

3) $\ti \nu$ is a 
 set    of isomorphisms $\{\nu_i: A_i\xrightarrow{\sim } A^{\vee}_{-i+a}\}_{i\in \Z/m\Z}$ 
 in $AV_{(p)}$ which are compatible with $\alpha_i$, and is a $\Z_{(p)}^\times$-homogeneous polarization,

4)  $\eta$ is a prime-to-$p$ level structure of type $K^p$ on $\{(A_i, \alpha_i, \ti \nu)\}_{i}$.

 For $K^p\subset \GG(\AA^p_f)$ sufficiently small, 
 this functor is represented by a  scheme $\SSh_{KK^p}({\rm GSp}(\mathbb V, \psi), S^{\pm})$ over $\Z_p$. The limit over $K^p$ is easily seen to satisfy (a) while (b) follows from the N\'eron-Ogg-Shafarevich criterion for good reduction of abelian varieties. Next, we discuss (c).

 For  a $\Z_p$-algebra $R$, set $\Lambda_{i, R}:=\Lambda_i\otimes_{\Z_p}R$ 
 and denote by $a_{i, R}: \Lambda_{i, R}\to \Lambda_{i+1, R}$ the map induced by the inclusion $\Lambda_i\to\Lambda_{i+1}$ and by $(\ ,\ )_i: \Lambda_{i, R}\times \Lambda_{-i+a, R}\to R$ the perfect $R$-bilinear pairing 
induced by $\psi$ since $\Lambda^\vee_i=\Lambda_{-i+a}$.  

  A result of G\"ortz \cite{GortzSymplectic} implies that, in this case, the local model ${\rm M}_K(G, \{\mu\})$ 
 of \cite{PaZhu} represents the functor which sends a $\Z_p$-algebra $R$ to the set of sequences $\{\F_i\}_{i\in \Z}$, where:
 $\F_i\subset \Lambda_{i, R}$ is an $R$-submodule which is Zariski locally
 a direct summand of $\Lambda_{i, R}$ of rank $g$, and, for all $i$, we have:  
 $ a_i(\F_i)\subset \F_{i+1}$,   $(\F_i, \F_{-i+a})_i=0$, and   $\F_{i-m}=p\F_i\subset p\Lambda_{i, R}=\Lambda_{i-m, R}$.

The morphism $\lambda$ is defined by first showing that the ``crystalline" system of modules, homomorphisms and pairings, obtained by applying
the Dieudonn\'e functor to $\{(A_i, \alpha_i, \nu_i)\}$, is locally isomorphic to the corresponding 
``Betti" system
given by the lattice chain $\{\Lambda_i\}$ and  $\psi$. Then,  $\lambda$ sends $\{(A_i, \alpha_i, \ti\nu, \eta)\}$
to $\F_i$ given by the Hodge filtration of $A_i$. It is smooth since,  by 
Grothendieck-Messing theory, deformations of an abelian variety are determined by lifts of its Hodge filtration. This argument first appeared in \cite{deJongGamma0} and \cite{DelignePappas}. It extends to most PEL type cases
\cite{RapZinkBook}. In general, we need a different approach which we explain next.

\subsection{} 
We now discuss general results.
Assume that $p$ is odd and that $(T)$ holds.
Assume also  that the local models ${\rm M}_K(G, \{\mu\})$ 
are as defined in \cite{PaZhu}.  The next result  is shown
in \cite{KisinPappas}, following earlier work of Kisin \cite{KisinJAMS}.

 \begin{ithm}\label{thm1}
 (i) Assume   that $(\GG, X)$ is of abelian type.
Then  Conjecture \ref{conjintegral}  with (c) replaced by (c'), is true 
for $(\GG, X)$ and $K$.

 (ii) Assume   that $(\GG, X)$ is of Hodge type and that the parahoric subgroup $K$ is equal to the stabilizer of a point 
 in the Bruhat-Tits building of $G$. 
Then   Conjecture \ref{conjintegral} is true for $(\GG, X)$ and $K$.
 \end{ithm}

 \subsection{}\label{proofs} 
 Let us try to explain some ideas in the proofs. 
 
 We will first discuss (ii). So,  the Shimura datum $(\GG, X)$ is of Hodge type and the parahoric subgroup $K\subset G(\Q_p)$ is the stabilizer of a point $z$ in the building $\B(G, \Q_p)$, \emph{i.e.} the stabilizer group scheme  is already connected,
$\Gg=\Gg_z=\Gg_z^\circ$. 

After adjusting the symplectic representation $\rho: \GG\hookrightarrow {\rm GSp}(\mathbb V, \psi)$ which gives the Hodge embedding, we can find a $\Z_p$-lattice $\Lambda\subset V=\mathbb V_{\Q_p} $ on which $\psi$ takes $\Z_p$-integral values, such that
$\rho$ induces
 
1) a closed group scheme embedding
$
\rho_{\Z_p}: \Gg=\Gg_z\hookrightarrow {\rm GL}(\Lambda),
$ and

2) an equivariant closed embedding 
$
\rho_*: {\rm M}_K(G, \{\mu\})\hookrightarrow {\rm Gr}(g, \Lambda^\vee)\otimes_{\Z_p}\O_E.
$

Here  ${\rm Gr}(g, \Lambda^\vee)$ is the Grassmannian over $\Z_p$
with $R$-points given by $R$-submodules $\F\subset \Lambda^\vee_R$ which are locally direct summands  of rank $g=\dim_\Q(\mathbb V)/2$. Finding $\Lambda$ is subtle and uses  that the representation $\rho$ is \emph{minuscule}.

 Choose a $\Z$-lattice $\VV_\Z\subset \VV$ such that $\VV_{\Z_{(p)}}:=\VV_\Z\otimes_\Z\Z_{(p)}=\VV\cap \Lambda$.
We can now find $\{s_{\alpha}\}_\alpha\subset \VV_{\Z_{(p)}}^\otimes\subset \VV^\otimes$ such that the functor
\[
 R\mapsto \GG_{\Z_{(p)}}(R)=\{g\in {\rm GL}(\VV_{\Z_{(p)}}\otimes_{\Z_{(p)}} R)\ |\ g\cdot s_\alpha=s_\alpha, \ \hbox{\rm for all $\alpha$}\}
\]
gives the unique flat group scheme $\GG_{\Z_{(p)}}$ over $\Z_{(p)}$ that extends both $\GG$ and $\Gg$. 
Here, $\VV^\otimes:=\bigoplus_{r, s\geq 0}( V^{\otimes r}\otimes (V^\vee)^{\otimes s})$, similarly for $\VV_{\Z_{(p)}}^\otimes$.

Now fix $K^p\subset \GG(\AA^p_f)$ small enough so, in particular, $V_{\hat\Z}:=V_\Z\otimes_\Z\hat\Z$ is stable by the action
of $K^p$; then $V_{\hat\Z}$ is $\eK=KK^p$-stable also. We can find $K'^p\subset {\rm GSp}(\VV)(\AA^p_f)$
such that $\rho$ gives a closed embedding
\[
\iota: \Sh_{KK^p}(\GG, X)\hookrightarrow \Sh_{K'_pK'^p}({\rm GSp}(\VV), S^{\pm})\otimes_{\Q}\mathbb E,
\]
where $K'_p$ is the subgroup of ${\rm GSp}(V, \psi)$ that stabilizes the lattice $\Lambda$, and such that
$\VV_{\hat\Z}$ is also
stable by the action of $\eK'=K'_pK'^p$.

We define $\SSh_{KK^p}(\GG, X)$ to be  the normalization of the
reduced Zariski closure of the image of $\iota_E$ in the integral model $\SSh_{K'_pK'^p}({\rm GSp}(\VV), S^{\pm})\otimes_{\Z_p}\mathbb \O_{E}$ (which is a moduli scheme, as in \S \ref{exampleSympl2}) and set 
$\SSh_{K}(\GG, X):=\varprojlim_{K^p}\SSh_{KK^p}(\GG, X)$. 
Checking property (a) is straightforward and (b) follows from the N\'eron-Ogg-Shafarevich criterion for good reduction. The hard work is in showing (c).

Since $\VV_{\hat\Z}$ is stable by the action of $\eK'$, which is sufficiently small, we have a universal abelian scheme  over $\Sh_{\eK'}({\rm GSp}(\VV), S^{\pm})$ which we can restrict via $\iota$ to obtain 
 an abelian scheme $\A$ over 
$\Sh_{\eK}(\GG, X)$. By construction, the tensors $s_\alpha$ above give sections $s_{\alpha, \rm B}$ of the local system  of $\Q$-vector spaces over $\Sh_\eK(\GG, X)(\C)$ with fibers ${\rm H}^1_{\rm B}(\A_v(\C), \Q)^\otimes$ given using 
 the first  Betti (singular) cohomology of the fibers $\A_v(\C)$ of the
 universal abelian variety. These are ``Hodge cycles" in the sense that they are of type $(0,0)$ for
 the induced Hodge structure on ${\rm H}^1_{\rm B}(\A(\C), \Q)^\otimes$.
 We are going to chase these around using various comparison isomorphisms.

Let $\kappa\supset \mathbb E$ be a field with  an embedding  $\sigma:\bar\kappa\hookrightarrow \C$ over $\mathbb E$
of its algebraic closure. Suppose $x\in \Sh_\eK(\GG, X)(\kappa)$ and let $\A_x$ be the corresponding abelian variety over $\kappa$. There are comparison isomorphisms
$
{\rm H}^1_B(\A_x(\C), \Q)\otimes_\Q\C\simeq {\rm H}^1_{\rm dR}(\A_x)\otimes_{\kappa,\sigma} \C$, and
$
{\rm H}^1_B(\A_x(\C), \Q)\otimes_\Q\Q_l\simeq {\rm H}^1_{\rm et}(\A_{x}\otimes_\kappa\bar\kappa, \Q_l),
$
 for any prime $l$.
Set $s_{\alpha,  B, x}$ be the fiber of $s_{\alpha, B}$ over the complex point $\sigma(x)$
and denote by $s_{\alpha, {\rm dR}, x}$, $s_{\alpha, {\rm et}, x}$, the images of 
$s_{\alpha, B, x}$ under these isomorphisms. The tensors $s_{\alpha, {\rm et}, x}$ are independent of the embedding $\sigma$, are fixed by the action of ${\rm Gal}(\bar\kappa/\kappa)$, and it follows from Deligne's ``Hodge implies absolute Hodge" theorem that 
$s_{\alpha, {\rm dR}, x}$ lie in ${\rm H}^1_{\rm dR}(\A_x)^\otimes$ and, in fact,
in $F^0({\rm H}^1_{\rm dR}(\A_x)^\otimes)$ (\hspace{-0.001cm}\cite{KisinJAMS}). 

For $l=p$, the corresponding ${\rm Gal}(\bar\kappa/\kappa)$-invariant tensors (``Tate cycles") $s_{\alpha, {\rm et}, x}$ lie in the lattice
${\rm H}^1_{\rm et}(\A_x\otimes_\kappa\bar\kappa, \Z_p)^\otimes=T_p(\A_{x})^\otimes$, 
where $T_p(\A_{x})$ is the 
$p$-adic Tate module $\varprojlim_n (\A_x\otimes_\kappa\bar\kappa)[p^n]$.
In order to control the local structure
of $\SSh_{KK^p}(\GG, X)$ and eventually relate it to ${\rm M}_K(G, \{\mu\})$, we need to employ some form of crystalline deformation theory and so we also have to understand the ``crystalline realization" of our tensors. 
Assume that $\kappa=F$ is a finite extension of $E$ and that the abelian variety $\A_x$ has good reduction. Then we obtain an $\O_F$-valued
point $\ti x$ of $\SSh_{KK^p}(\GG, X)$ that extends $x$ and reduces to $\bar x$. 
  A key point is to show that the Tate cycles $s_{\alpha, {\rm et}, x}$ for $l=p$
give, via the \'etale/crystalline comparison, corresponding ``nice" \emph{integral} crystalline cycles $s_{\alpha, {\rm cris}, \bar x}\in \DD(\A_{\bar x})^\otimes$ on the Dieudonn\'e
module $\DD(\A_{\bar x})$ of the abelian variety $\A_{\bar x}$ over the residue field of $F$. 
For this we need a suitably functorial construction that relates crystalline $p$-adic Galois representations 
of ${\rm Gal}(\bar F/F)$ to Frobenius semilinear objects integrally, and without restriction
on the absolute ramification of $F$. This is provided
by the theory of Breuil-Kisin modules \cite{KisinFF, KisinJAMS}.

Let $F_0$ be the maximal unramified extension of $\Q_p$ contained in $F$. Denote by $W_0$ the ring of integers of $F_0$ and  by $\phi:\fS\to \fS$ the lift of Frobenius
with $\phi(u)=u^p$. Choose a uniformizer $\pi$ of $F$ which is a root of an Eisenstein polynomial $E(u)\in W_0[u]$. 
A Breuil-Kisin module $(\frakM, \Phi)$ is a finite free $\fS$-module $\frakM$
 with an isomorphism $\Phi: \phi^*(\frakM)[1/E(u)]\xrightarrow{\sim} \frakM[1/E(u)]$, 
where $\phi^*(\frakM):=\fS\otimes_{\phi,\fS}\frakM$.
Kisin has constructed a fully faithful tensor functor 
$
T\mapsto \frakM(T)
$ from the category of  ${\rm Gal}(\bar F/F)$-stable $\Z_p$-lattices in crystalline $\Q_p$-representations to the category of Breuil-Kisin modules.  

Now let $T^\vee={\rm Hom}_{\Z_p}(T,\Z_p)$ be the linear dual of the $p$-adic Tate module $T=T_p(\A_{x})$. Then
there are natural isomorphisms \cite{KisinFF}
\[
\phi^*(\frakM(T^\vee)/u\frakM(T^\vee))\simeq \DD(\A_{\bar x}),\quad \phi^*(\frakM(T^\vee))\otimes_{\fS, u\mapsto \pi}\O_F\simeq {\rm H}^1_{\rm dR}(\A_{\ti x}).
\]
Applying Kisin's functor to the Tate cycles $s_{\alpha, {\rm et}, x}$ gives Frobenius invariant tensors $\tilde s_\alpha\in \frakM(T_p(\A_{x})^\vee)^\otimes$.
By using the isomorphism above, these give the crystalline cycles $s_{\alpha, {\rm cris}, \bar x}\in \DD(\A_{\bar x})^\otimes$.
The main result now is: 

\begin{ithm}\label{keylemma}
There is a $\fS$-linear isomorphism 
\[
\beta: \Lambda^\vee\otimes_{\Z_p} \fS\xrightarrow{\sim} \frakM(T_p(\A_{x})^\vee)
\] such that $\beta^\otimes$ takes  $s_\alpha\otimes 1$ to $\tilde s_\alpha$, for all $\alpha$.
\end{ithm}

In addition to the properties of Kisin's functor, the important input in the proof is the
statement that all $\Gg$-torsors over $\Spec(\fS)-\{(0,0)\}$ are trivial. This uses crucially that $\Gg=\Gg_z^\circ$ is parahoric.

Using $\beta$ and the above, we obtain isomorphisms 
\[
\beta_{\rm dR}: \Lambda^\vee\otimes_{\Z_p} \O_F\xrightarrow{\sim} {\rm H}^1_{\rm dR}(\A_{\ti x}),\quad \beta_{\rm cris}: \Lambda^\vee\otimes_{\Z_p} W_0\xrightarrow{\sim} \DD(\A_{\bar x}),
\]
such that $ \beta^\otimes_{\rm dR}(s_\alpha\otimes 1)=s_{\alpha, {\rm dR}, \ti x}$,  resp. $ \beta^\otimes_{\rm cris}(s_\alpha\otimes 1)=s_{\alpha, {\rm cris}, \bar x}$\footnote{In the case that $K$ is hyperspecial
and $F=F_0$, the existence of such an isomorphism $\beta_{\rm cris}$ was conjectured by Milne 
and was shown by Kisin \cite{KisinJAMS} and Vasiu \cite{VasiuMotivic}.}.The inverse image
\[
\beta^{-1}_{\rm dR}(F^0(\A_{\ti x}))\subset \Lambda^\vee\otimes_{\Z_p} \O_F
\]
of the deRham filtration   $F^0(\A_{\ti x})={\rm H}^0(\A_{\ti x}, \Omega^1_{\A_{\ti x}/\O_F})\subset {\rm H}^1_{\rm dR}(\A_{\ti x})$ now gives an $\O_F$-point of the Grassmannian 
${\rm Gr}(g, \Lambda^\vee)$. This is easily seen to factor through $\rho_*$ to the local model ${\rm M}_K(G, \{\mu\})$
and defines the image of $\lambda$ on $\tilde x\in \SSh_{\eK}(\GG, X)(\O_F)$.

Next, we use  $\beta$ to construct a versal deformation of the abelian scheme $\A_{\tilde x}$ over the  completion $S=\widehat{\O}_{{\rm M}_K(G,\{\mu\}), \bar y}$ of the local ring of ${\rm M}_K(G,\{\mu\})$ at $\bar y=\lambda(\bar x)$. This deformation is equipped with Frobenius invariant crystalline tensors that appropriately extend $\tilde s_\alpha$. It is constructed using Zink's theory of displays \cite{Zinkdisplay}. In this, (connected) $p$-divisible groups over a $p$-adic ring $S$ are described by ``displays",  which are systems of Frobenius modules  over the ring of 
Witt vectors $W(S)$.
In our case, we give a display over $W(S)$ whose Frobenius semilinear maps 
are, roughly speaking, valued in the group scheme $\Gg$. Our construction is inspired by the idea that the local model ${\rm M}_K(G, \{\mu\})$ should also appear inside a Witt  affine Grassmannian  (see \S \ref{resultsLocal} and also \S \ref{displays} for a similar
construction).  

 Finally,
using a result of Blasius and Wintenberger and parallel transport for the Gauss-Manin connection, we show that there exists 
a  homomorphism $S\to \widehat\O_{\SSh_{\eK}(\GG, X),\bar x}$
which matches $\tilde x$ and $\tilde y$ and then has to be an isomorphism. Property (c') and then (c) follows.

To deal with the case that $(\GG, X)$ is of abelian type, we follow Deligne's strategy \cite{DeligneCorvallis} for constructing canonical models for these types. Deligne relates a connected component of the Shimura variety for $(\GG, X)$ to one of the Hodge type $(\GG_1, X_1)$ (as in the definition of abelian type above). In this, he uses an action of $\GG^{\rm ad}(\Q)\cap \GG^{\rm ad}(\R)^+$ on these varieties. The argument extends after
giving a moduli interpretation for this action. Then, $\SSh_{K}(\GG, X)$ is constructed, by first taking a quotient of
a connected component of $\SSh_{K_1}(\GG_1, X_1)$ by a finite group action 
(which is shown to be free), and then ``inducing" from that quotient.   

The strategy for these proofs is due to Kisin \cite{KisinJAMS} who showed the results when $K$ is hyperspecial (\emph{i.e.} $\Gg$ is reductive). Then ${\rm M}_K(G, \{\mu\})$ 
is the natural smooth model of the homogeneous space $X_\mu$ and $\SSh_{\eK}(\GG, X)$ is smooth over $\O_E$. (The condition on $\pi_1(G^{\rm der})$ is not needed then. Also, the case $p=2$ is treated by Kim and Madapusi Pera \cite{KimMadapusi}.) In this important hyperspecial case, the integral models $\SSh_{K}(\GG, X)$ can be shown to be 
``canonical", in the sense of Milne: They can be characterized as the unique, up to isomorphism, regular, formally smooth schemes that satisfy (a), (c'), and an extension property stronger than (b) in which $\Spec(R)$ is replaced by 
any regular, formally smooth scheme over $\O_E$ (see \cite{KisinJAMS}, \cite{VasiuZink}). Hence, they are also independent of the various choices 
made in their construction. It is unclear if schemes satisfying (a), (b) and (c) as in the conjecture, 
with ${\rm M}_K(G,\{\mu\})$ given a priori, are uniquely determined in the general parahoric case.
In \cite{KisinPappas}, we show this in a few cases.
Regardless, we conjecture that the schemes $\SSh_{K}(\GG, X)$ produced by the above results are
independent of the choices made in their construction. 

In the hyperspecial case, there is also earlier work of Vasiu \cite{VasiuAJM} who pursued a different approach (see also \cite{Moonen}). Vasiu applied directly an integral comparison homomorphism between $p$-adic \'etale cohomology and crystalline cohomology due to Faltings, to understand carefully chosen crystalline tensors of low degree that should be enough to control the group.

Let us note here that, in contrast to all the previous theory, there are few results when $K$ is a  subgroup smaller than parahoric, even in PEL cases:
When $K$ is the pro-$p$ radical of an Iwahori subgroup,  combining the above with
results of Oort-Tate or Raynaud on $p$-torsion group schemes, often leads to well-behaved models
(\hspace{-0.001cm}\cite{Pappas95}, \cite{HainesRapoport}). 
Also, when the deformation theory is controlled by a formal group 
(or ``$\O$-module") of dimension $1$ (as is in the case of modular curves) the notion of Drinfeld level structure can be used to
describe integral models for all level subgroups (\hspace{-0.001cm}\cite{KatzMazur}, \cite{HarrisTaylor}).

\section{\bf Reductions: Singularities and points. }

\subsection{}\label{KR} 
\  \emph{Singularities: Kottwitz-Rapoport stratification and nearby cycles.} 
To fix ideas, we assume that the assumptions of  Theorem \ref{thm1} (ii) are satisfied so that the result applies.
In particular, $(\GG, X)$ is of Hodge type. 
Set $\SSh:=\SSh_{\eK}(\GG, X)$, ${\rm M}:={\rm M}_K(G,\{\mu\})$, $S:=\Spec(\O_E)$ and $d:=\dim(\Sh_{\eK}(\GG, X))$.
 
Since the morphism $\lambda$ is smooth, Theorem \ref{normalThm} implies that special fiber $\SSh\otimes_{\O_E}k_E$ is reduced,
and that we obtain a stratification of  $\SSh\otimes_{\O_E}\BFf$ 
by locally closed smooth strata ${\SSh}_{[w]}=\lambda^{-1}_\eK(S_{[w]})$, for $[w]\in {\rm Adm}_{\ti K}(\{\mu\})$.
Furthermore:

\begin{ithm} 
 For each $[w]\in {\rm Adm}_{\ti K}(\{\mu\})$, the Zariski closure $\ov {\SSh}_{[w]}$   is equal to the union
$
\bigcup_{[w']\leq [w]} \SSh_{[w']}
$ and is normal and Cohen-Macaulay. 
\end{ithm}

Now pick a prime $l\neq p$ and consider the complex  of nearby cycles 
${\rm R}\Psi(\SSh, \bar \Q_l)$ of $\SSh\to S$ on $\SSh\otimes_{\O_E}\BFf$.
This is obtained from the nearby cycles 
${\rm R}\Psi({\rm M}, \bar \Q_l)$ of ${\rm M}\to S$ by pulling back along the smooth $\lambda$.
 Set $f=[k_E:\Ff_p]$ and $q=p^f$. 
For each $r\geq 1$, the semi-simple trace of Frobenius \cite{RapoportAnnArbor} defines a function \vspace*{-0.15cm}
\[
{\psi}_{\eK, r}: \SSh({\mathbb F}_{q^r})\to \bar\Q_l;\quad \psi_{\eK, r}(\bar x)={\rm Tr}^{\rm ss}({\rm Frob}_{\bar x}, {\rm R}\Psi(\SSh, \bar\Q_l)_{\bar x}).\vspace*{-0.15cm}
\]
The smoothness of $\lambda$ and $\Gg$, and Lang's lemma, gives that ${\psi}_{\eK, r}$ factors through 
\vspace*{-0.15cm}
\[
\lambda_{\eK}({\mathbb F}_{q^r}) : \SSh({\mathbb F}_{q^r})\to \Gg(\Ff_{q^r})\backslash {\rm M}_K(G, \{\mu\})(\Ff_{q^r}). 
\]

Denote by $\Q_{p^n}$ the unramified extension of $\Q_p$ of degree $n$ contained in $L$; let $\Z_{p^{n}}\simeq W(\Ff_{p^{n}})$ be its integers. Set $K_r=\Gg(\Z_{q^{r}})$. If $G$ is quasi-split over $\Q_{q^r}$,
then we have a bijection
$
\Gg(\Ff_{q^r})\backslash {\rm M}_K(G, \{\mu\})(\Ff_{q^r})\cong K_r\backslash G(\Q_{q^r})/K_r,
$
and so the map $\lambda_{\eK}({\mathbb F}_{q^r})$ gives
$
\lambda_{\eK, r} : \SSh({\mathbb F}_{q^r})\to K_r\backslash G(\Q_{q^r})/K_r .
$
\smallskip

By combining  results of \cite{PaZhu} on ${\rm R}\Psi({\rm M}, \bar \Q_l)$ with the above,  we obtain:

\begin{ithm}
a) Suppose that $G$ splits over the finite extension 
$F/\Q_p$. Then the inertia subgroup $I_{EF}\subset {\rm Gal}(\bar \Q_p/EF)$ acts unipotently on ${\rm R}\Psi(\SSh, \bar\Q_l)$.

  b) Suppose that $G$ is quasi-split over $\Q_{q^r}$. 
Then   the semi-simple trace of Frobenius ${\psi}_{r, \eK}: \SSh({\mathbb F}_{q^r})\to \bar\Q_l$ 
is a composition \vspace*{-0.1cm}
\[
 \SSh({\mathbb F}_{q^r})\xrightarrow{\lambda_{\eK, r} } K_r\backslash G(\Q_{q^r})/K_r \xrightarrow{z_{\mu, r} } \bar\Q_l
\]
 where $z_{\mu, r}$ is in the center of the Hecke algebra
 $\bar\Q_l[K_r\backslash G(\Q_{q^r})/K_r] $ with multiplication by convolution.
 \vspace{0.05cm}

c)  Suppose that $G$ is split over $\Q_{q^r}$.  Then in (b) above, we can take $z_{\mu, r}$
such that $q^{-d/2}z_{\mu, r}$ is a Bernstein function for $\mu$. 
\end{ithm}

 Part (c) has been conjectured by Kottwitz,
see \cite[(10.3)]{RapoportGuide} and also work of Haines \cite{HainesBernsteinCenter} 
for more details and an extension of the conjecture. Let us mention  that the study of ${\rm R}\Psi({\rm M}, \bar \Q_l)$
in \cite{PaZhu} (as also the proof of Theorem \ref{normalThm}), uses techniques from the 
theory of the geometric Langlands correspondence. The ``Hecke central" statement was first shown by
Gaitsgory 
 for split groups over function fields \cite{GaitsgoryInv}, and by Haines and Ng\^o for unramified unitary and symplectic groups \cite{HainesNgoNearby}. Gaitsgory's result was, in turn, inspired by Kottwitz's conjecture.

\subsection{}
\emph{Points modulo $p$.}  Under certain assumptions, Langlands and Rapoport  gave a conjectural description of the set of $\BFf$-points of a model of a Shimura variety together with its actions by Frobenius and Hecke operators \cite{LanglandsRapoport}. The idea, very roughly, is as follows (see also \cite{Langlands}, \cite{KottJAMS}). Since, as suggested by Deligne, most Shimura varieties are supposed to be moduli spaces of certain motives, we should be describing this set via representations of the fundamental groupoid attached to the Tannakian category of motives over $\BFf$. A groupoid $\frakQ$ which should be this fundamental groupoid (and almost is, assuming the Tate conjecture and other standard conjectures) can be constructed explicitly. Then ``$\GG$-pseudo-motives" are given by groupoid homomorphisms
$\phi: \mathfrak Q\to \GG$. The choice of the domain $X$ imposes restrictions,
and we should be considering only $\phi$  which are ``admissible". Then, we also give $p$ and prime-to-$p$ level structures on these $\phi$,  as we do for abelian varieties. 

The Langlands-Rapoport conjecture  was corrected, modified and extended along the way (\hspace{-0.001cm}\cite{Reimann},\cite{KisinLR}). The conjecture makes 
better sense when it refers to a specific integral model of the Shimura variety $\Sh_K(\GG, X)$.
When $K$ is hyperspecial, $p$ odd, and $(\GG, X)$ of abelian type, an extended version of the conjecture was essentially proven (with a caveat, see below) by Kisin \cite{KisinLR}, for the canonical integral models 
constructed in \cite{KisinJAMS}. In fact,  the conjecture also makes sense when $K$ is parahoric, for the integral models of \cite{KisinPappas}.   More precisely, let us suppose  that the assumptions of Theorem \ref{thm1} are satisfied and that $\SSh_K(\GG, X)$
is provided by the construction in the proof. 
\begin{iconjecture}\label{LRconj}
There is a $\langle\Phi_E\rangle\times \GG(\AA_f^p)$-equivariant bijection
\begin{equation*}\label{LReq}
\SSh_K(\GG, X)(\BFf)\xrightarrow{\sim} \bigsqcup\nolimits_{[\phi]} \varprojlim\nolimits_{K^p} I_\phi(\Q)\backslash (X_p(\phi)\times (X^p(\phi)/K^p))\leqno(LR)
\end{equation*}
where $\Phi_E$ is the Frobenius over 
$k_E$ and the rest of the notations and set-up follow \cite[\S 9]{RapoportGuide}, \cite[(3.3)]{KisinLR} (taking into account the remark (3.3.9) there).
\end{iconjecture} In particular: 
The disjoint union is over a set of equivalence classes $[\phi]$ of admissible $\phi:\frakQ\to \GG$,
where one uses Kisin's definition of ``admissible'' for non simply connected derived group. 
(The classes $[\phi]$ often correspond to isogeny classes of abelian varieties with additional structures.) 
We have $I_\phi={\rm Aut}(\phi)$, an algebraic group over $\Q$. The set $X^p(\phi)$ is a right $\GG(\AA_f^p)$-torsor
and the set $X_p(\phi)$ can be identified with the affine Deligne-Lusztig set $X_K(\{\mu^{-1}\}, b)$ for  $b\in G(L)$ obtained
from $\phi$. 
The group $I_\phi(\Q)$ acts on $X^p(\phi)$ on the left, and there is an injection $
I_\phi(\Q)\to J_b(\Q_p)$ which also produces a left action on $X_p(\phi)=X_K(\{\mu^{-1}\}, b)$; the quotient is by the diagonal action.

\begin{ithm} (Kisin \cite{KisinLR})
Assume that $K$ is hyperspecial, $p$ odd, and $(\GG, X)$ of abelian type. Then there is a
bijection as in (LR) respecting the action of $\langle\Phi_E\rangle\times\GG(\AA_f^p)$ on both sides, 
but with the action of 
$I_\phi(\Q)$ on $X_p(\phi)\times (X^p(\phi)/K^p)$ obtained from the natural diagonal action above
by conjugating by a (possibly trivial) element $\tau(\phi)\in I^{\rm ad}_\phi(\AA_f)$.
\end{ithm}
 
Due to lack of space, we will omit
an account of any of the beautiful and subtle arguments in Kisin's proof, or in related earlier work by Kottwitz (\hspace{-0.001cm}\cite{KottJAMS}) and others. 
We hope that the interested reader will consult the original papers. Let us just mention here that R. Zhou \cite{Zhou} has recently made some progress towards the proof of Conjecture \ref{LRconj}.

\section{\bf Local theory: Formal schemes.} \label{localShi}

\subsection{} 
We now return to the local set-up as in \S \ref{Local}.  
Let $(G, \{\mu\})$ be a local Shimura pair
and $[b]\in B(G)$ a $\sigma$-conjugacy class. Fix a parahoric subgroup $K\subset G(\Q_p)$
and assume that the local model ${\rm M}_K(G,\{\mu^{-1}\})$ is defined. 
 
\begin{iconjecture}\label{conjformal} There exists a formal scheme $\frakX_K(\{\mu\}, b)$  over ${\rm Spf}(\O_E)$ with $J_b(\Q_p)$-action, which is locally formally of finite type, and is such that:

a) There is a $\langle \Phi_E\rangle \times J_b(\Q_p)$-equivariant bijection between $\frakX_K(\{\mu\}, b)(\BFf)$ 
and the affine Deligne-Lusztig set $X_K(\{\mu\}, b)$.

b) For any $v\in X_K(\{\mu\}, b)$, there is $w\in {\rm M}_K(G,\{\mu^{-1}\})(\BFf)$ such that the completion of $\frakX_K(\{\mu\}, b)$ at the point corresponding to $v$ via (a) is isomorphic to  the completion of  ${\rm M}_K(G,\{\mu^{-1}\})$ at $w$. 
\end{iconjecture}

Assuming $\frakX_K(\{\mu\}, b)$ exists, let $\mathscr X_K(\{\mu\}, b)$ be its underlying reduced scheme which is locally of finite type over $\Spec(k_E)$.
(Then  $\mathscr X_K(\{\mu\}, b)$ could be called an \emph{affine Deligne-Lusztig ``variety''}.)
 Its ``perfection" $\mathscr X_K(\{\mu\}, b)^{\rm perf}$, has been constructed by Zhu \cite{ZhuWitt}, using a Witt vector affine flag variety.
 The rigid-analytic fiber $\frakX_K(\{\mu\}, b)^{\rm rig}$ over $E$ would be the \emph{local Shimura variety} for $(G,\{\mu\}, [b])$ and level $K$ whose existence is expected by Rapoport and Viehmann \cite{RapoportViehmann}. 

Let us take $G={\rm GL}_n$, $\mu_d(z)={\rm diag}(z,\ldots ,z, 1, \ldots ,1)^{-1}$, with $d$ copies of $z$,  and $K={\rm GL}_n(\Z_p)$.  For simplicity,  set $W=W(\BFf)$. The formal schemes $\frakX_{\GL_n}:=\frakX_K(\{\mu_d\}, b)$ were constructed by Rapoport and Zink \cite{RapZinkBook}: Fix a $p$-divisible group $H_0$ of height $n$ and dimension $d$ over $\BFf$ with  Frobenius $F=b\cdot\sigma$ on the rational Dieudonn\'e module. Then the formal scheme
$\frakX_{\GL_n}\hat\otimes_{\Z_p}W$ represents (\hspace{-0.001cm}\cite[Theorem 2.16]{RapZinkBook}) the functor which sends a $W$-algebra $R$ with $p$ nilpotent on $R$, to the set of isomorphism classes
of pairs $(H, \tau)$, where:   $H$ is a $p$-divisible group over $R$, and 
  $\tau: H_0\otimes_k R/pR\dashrightarrow H\otimes_R R/pR$ is a quasi-isogeny. 
The formal scheme $\frakX_{\GL_n}$ over ${\rm Spf}(\Z_p)$ is obtained from this by descent. 
This construction of $\frakX_K(\{\mu\}, b)$ generalizes  to $(G,\{\mu\})$
that are of ``EL" or ``PEL" type, for many parahoric $K$; these types are defined similarly to  (and often arise from) global Shimura data $(\GG, X)$ of PEL type.  
 
 \subsection{}
The integral models $\SSh_\eK(\GG, X)$ of Conjecture \ref{conjintegral} and the formal schemes  $\frakX_K(\{\mu\}, b)$ above, should be intertwined via the bijection $(LR)$ 
 as follows: 

Suppose that $(\GG, X)$ is a Shimura datum which produces the local Shimura pair $(G, \{\mu^{-1}\})$. Take $\bar x\in \SSh_K(\GG, X)(\BFf)$ with corresponding $b\in G(L)$. Then, there should be a morphism of ${\rm Spf}(\O_E)$--formal schemes
\vspace{-0.1cm}
\[
i_{\bar x}:  \frakX_K(\{\mu\}, b)\times \GG(\AA_f^p)/K^p\xrightarrow{ } \widehat{\SSh}_\eK(\GG, X):=\varprojlim\nolimits_{n}\SSh_\eK(\GG, X)\otimes_{\O_E}\O_E/(p^n), \vspace{-0.1cm}
\]
 which, on $\BFf$-points is given by $(LR)$ and
surjects on the ``isogeny class"
$[\phi_0]$ of $\bar x$, and induces isomorphisms on the formal completions at closed points.

 An interesting special case is when the $\sigma$-conjugacy class $[b]$ 
is basic (\hspace{-0.001cm}\cite{KottIsocrystalsII}). Then the image $Z$ of $i_{\bar x}$ should be closed
in the special fiber of $\SSh_\eK(\GG, X)$  and $i_{\bar x}$ 
should be giving a ``$p$-adic uniformization" (\emph{cf.} \cite{RapZinkBook})
\vspace{-0.1cm}
\[
 \ I_{\phi_0}(\Q)\backslash \frakX_K(\{\mu\}, b)\times (\GG(\AA_f^p)/K^p)\xrightarrow{\sim } \widehat{\SSh}_\eK(\GG, X)_{/Z}, \leqno(U) \vspace{-0.1cm}
\]
of the completion of $\SSh_\eK(\GG, X)$ along $Z$.  

Assume that $(\GG, X)$ is of Hodge type and that the rest of the assumptions of  Theorem \ref{thm1} (ii) are also satisfied. Then we can hope to construct $\frakX_K(\{\mu\}, b)$
using the model $\SSh_\eK(\GG, X)$ as follows: 
First consider the fiber product
 \begin{equation*}
\begin{matrix}  {\mathfrak F}& \to & (\frakX_{{\rm GSp}}\times 1)\hat\otimes_{\Z_p}\O_E \\
\downarrow && \ \downarrow{i_{\bar x}}\\
\widehat{\SSh}_\eK(\GG, X) &\xrightarrow{\iota } & \widehat{\SSh}_{\eK'}({\rm GSp}(\VV), S^\pm)\hat\otimes_{\Z_p}\O_E,
\end{matrix}\vspace{-0.08cm}
 \end{equation*}
 where $\iota$ is an appropriate Hodge embedding as in \cite{KisinJAMS}, or \cite{KisinPappas},
 and $\frakX_{\GSp}$ is the Rapoport-Zink formal scheme for the symplectic PEL type. 
 Then, provided we have the set map $X_K(\{\mu\}, b)\to \SSh_\eK(\GG, X)(\BFf)$
 sending $1\cdot\ti K$ to $\bar x$ (as predicted by the $(LR)$ conjecture), we can define $\frakX_K(\{\mu\}, b)$ to be the formal completion of $\frakF$ along the (closed) subset given by $X_K(\{\mu\}, b)$. The existence of 
 the ``uniformization" morphism $i_{\bar x}$ follows immediately.  Howard and the author applied  this idea to show:
 
 \begin{ithm}\label{HP} (\hspace{-0.001cm}\cite{HP2})
 Assume $p\neq 2$, $(\GG, X)$ is of Hodge type and $K$ is hyperspecial. 
 Choose a Hodge embedding $\rho$, a lattice $\Lambda$,
 and tensors $\{s_\alpha\}_\alpha$, as in the proof of Theorem \ref{thm1} (ii). 
 Suppose $\bar x\in \SSh_{\eK}(G, X)(\BFf)$ and $b$ are as above.  Then, $\frakX_K(\{\mu\}, b)$ satisfying Conjecture \ref{conjformal} can be defined
 as above, and only depends, up to isomorphism, on $(G, \{\mu\}, b, K)$ and 
 $\rho_{\Z_p}:\Gg\to \GL(\Lambda)$. 
 \end{ithm}
 
  For basic $[b]$, the uniformization $(U)$ also follows.
In the above result, the formal scheme
 $\frakX_K(\{\mu\}, b)\hat\otimes_{\O_E} W$ represents a functor on the category of 
 formally finitely generated adic $W$-algebras $(A, I)$ which are formally smooth over $W/(p^m)$, for some $m\geq 1$ as follows. Set $H_0$ for the $p$-divisible group $\A_{\bar x}[p^\infty]$. We send $(A, I)$ to the set of isomorphism classes of  triples $(H, \tau, \{t_\alpha\}_\alpha)$, where $H$ is a $p$-divisible group over $A$, $\tau: H_0\otimes_k A/(p, I)\dashrightarrow H\otimes_A A/(p, I)$ is a quasi-isogeny, and $t_\alpha$ are (Frobenius invariant) sections of the Dieudonn\'e crystal $\DD(H)^\otimes$. We require that the pull-backs $\tau^*(t_\alpha)$ agree in
 $\DD(H_0\otimes_k A/(p, I))^\otimes[1/p]$
 with the constant sections given by $s_{\alpha, {\rm crys}, \bar x}$,
 and also some more properties listed in \cite[(2.3.3), (2.3.6)]{HP2}.
For $K$ hyperspecial, a different, local, construction of the formal scheme $\frakX_K(\{\mu\}, b)$
(still only for ``local Hodge types") first appeared in work of W. Kim \cite{KimHodge}.

  \subsection{}\label{displays}
  In some cases, we can also give $\frakX_K(\{\mu\}, b)$ directly,
  by using a group theoretic version of Zink's displays.
  Assume that $K$ is hyperspecial, so  $\Gg$ is reductive over $\Z_p$. Set $\Gg_W=\Gg\otimes_{\Z_p}W$ and pick $\mu: \Gm_W\to \Gg_W$ in our conjugacy class. Set $L^+\Gg$, resp. $L\Gg$, for the ``positive Witt loop" group scheme, resp. ``Witt loop" ind-group scheme, representing $R\mapsto \Gg(W(R))$, resp. $R\mapsto \Gg(W(R)[1/p])$. 
   Let $\calH_\mu$ be the subgroup scheme of $L^+\Gg_W $ with $R$-points   given  
by $g\in \Gg(W(R))$ whose projection $g_0\in \Gg(R)$ lands in the $R$-points of the parabolic subgroup scheme $\calP_\mu\subset \Gg_W$ associated to $\mu $. 
We can define a   homomorphism
$
\Phi_{G, \mu}: \calH_\mu \to L^+\Gg_W
$
such that $\Phi_{G, \mu}(h)={ F}\cdot (\mu(p)\cdot h\cdot \mu(p)^{-1})$ in $\Gg(W(R)[1/p])$, with
$F$ given by the Frobenius $W(R)\to W(R)$ (\hspace{-0.001cm}\cite{OP}).
Let $R$ be a $p$-nilpotent $W$-algebra. 
Set 
\vspace*{-0.1cm}
 \[
\frakX'_\Gg(R)= \{(U, g)\in L^+\Gg(R)\times L\Gg(R)\ |\  g^{-1}bF(g)\buildrel{(*)}\over{=}U\mu^\sigma(p)\}/ \calH_\mu(R),
 \]
where $(*)$ is taken in $L\Gg(R)$ and the quotient is for the action given by
\[\vspace*{-0.1cm}
(U, g)\cdot h=(h^{-1} U  \Phi_{G,\mu}(h), g h).
\]

Denote by $\frakX_\Gg$ the \'etale sheaf associated to  $R\mapsto \frakX'_\Gg(R)$. 

If $R$ is perfect, then $F$ is an isomorphism and $W(R)$ is $p$-torsion free, so a pair $(U, g)$ with $g^{-1}bF(g)=U\mu^\sigma(p)$ is determined by $g$. Then
\[
\frakX'_\Gg(R)=\{g\in L\Gg(R)/L^+\Gg(R)\ |\ g^{-1}bF(g)\in L^+\Gg(R)\mu^\sigma(p)L^+\Gg(R)\}.
\]
In particular, $\frakX_\Gg(\BFf)=\frakX'_\Gg(\BFf)\cong  X_K(\{\mu\}, b)$ and the perfection of $\frakX_\Gg$ 
agrees with the space $\mathscr X_K(\{\mu\}, b)^{\rm perf}$ of Zhu \cite{ZhuWitt}. 
If  $\frakX_\Gg$   is representable by a formal
scheme, then this satisfies Conjecture \ref{conjformal} and gives
$\frakX_K(\{\mu\}, b)\hat\otimes_{\O_E}W$.
 
\begin{ithm}(\hspace{-0.001cm}\cite{OP})
Let $\rho_{\Z_p}: \Gg\hookrightarrow \GL_n$ be a closed group scheme
embedding. Suppose $\rho_{\Z_p}\circ \mu$ is minuscule and $\rho_{\Z_p}(b)$ has no zero slopes.
Then the restriction of $\frakX_\Gg$ to Noetherian $p$-nilpotent $W$-algebras is represented
by a formal scheme which is formally smooth and locally formally of finite type over ${\rm Spf}(W)$. 
\end{ithm}
 
The proof also gives that, if $\frakX_K(\{\mu\}, b)$ is defined as in Theorem \ref{HP} and 
$\rho_{\Z_p}(b)$ has no zero slopes, then 
$
\frakX_\Gg(R) =\frakX_K(\{\mu\}, b)(R),
$
for $R$  Noetherian.
Hence, then $\frakX_K(\{\mu\}, b)$ of Theorem \ref{HP} is independent of $\rho_{\Z_p}$.

 \providecommand{\bysame}{\leavevmode\hbox to3em{\hrulefill}\thinspace}
\providecommand{\href}[2]{#2}

{\footnotesize DEPT. OF MATH., MICHIGAN STATE UNIVERSITY, E. LANSING, MI 48824, USA}

\end{document}